\documentclass[10 pt, conference]{ieeeconf}  % Comment this line out

\IEEEoverridecommandlockouts                              % This command is only
\title{\LARGE \bf
Differentially Private Distributed Optimal Power Flow
}

\author{Vladimir Dvorkin$^{1}$, Pascal Van Hentenryck$^{2}$, Jalal Kazempour$^{1}$, Pierre Pinson$^{1}$% <-this % stops a space
\thanks{$^{1}$V. Dvorkin, J. Kazempour. and P. Pinson are with the Department of Electrical Engineering,
        Technical University of Denmark, Lyngby, Denmark.
        Email: {\tt \{vladvo,seykaz,ppin\}@elektro.dtu.dk}}%
\thanks{$^{2}$Pascal Van Hentenryck is with the School of Industrial and Systems Engineering, Georgia Institute of Technology, Atlanta, USA.
        Email: {\tt pascal.vanhentenryck@isye.gatech.edu}}%
}

\usepackage[english]{babel}
\usepackage[dvipsnames]{xcolor}
\usepackage{thmtools}
\usepackage{graphicx}
\usepackage[normalem]{ulem}
\usepackage{amsmath}
\usepackage{amssymb}
\usepackage{amsfonts}
\usepackage{enumerate}
\usepackage{euscript}
\usepackage{empheq}
\usepackage{mathtools}
\usepackage{tikz}
\usepackage{pgfplots}
\usepackage{euscript}
\usepackage{hyperref}
\usepackage{balance}
\usepackage{cite}
\usepackage{booktabs}
\usepackage{xfrac}
\usepackage{algorithm}
\usepackage{multirow}
\usepackage[noend]{algpseudocode}
\usetikzlibrary{arrows}
\usetikzlibrary{shapes}
\usetikzlibrary{patterns}
\usepgfplotslibrary{fillbetween}
\usetikzlibrary{arrows.meta}
\usepgfplotslibrary{colorbrewer}

\newcommand{\minimize}[1]{\underset{{#1}}{\text{min}}}
\newcommand{\argmin}[1]{\underset{{#1}}{\text{argmin}}}
\newcommand{\maximize}[1]{\underset{{#1}}{\text{max}}}

\newcommand{\braceit}[1]{\left({#1}\right)}
\newcommand{\st}{\text{s.t.}}
\newcommand{\set}{\EuScript}
\usepackage{mathtools}
\DeclarePairedDelimiterX{\norm}[1]{\lVert}{\rVert}{#1}
\newcommand*\circled[1]{\tikz[baseline=(char.base)]{
            \node[shape=circle,draw,inner sep=0.5pt] (char) {\footnotesize{#1}};}}
\newcommand\mydots{\hbox to 0.75em{.\hss.\hss.}}

\tikzset{myarr/.style={{Circle[black,length=4pt]}-{Circle[black,length=4pt]},shorten <=-2.5pt,shorten >=-2.5pt}}

\DeclareMathSymbol{\shortminus}{\mathbin}{AMSa}{"39}

\newcommand{\kplus}{k\scalebox{.7}{+}1}
\newcommand{\tplus}{t\scalebox{.7}{+}1}

\definecolor{red}{HTML}{e05a87} 
\definecolor{blue}{HTML}{035096}

\newtheorem{assumption}{Assumption}
\newtheorem{definition}{Definition}
\newtheorem{proposition}{Proposition}
\newtheorem{theorem}{Theorem}
\newtheorem{remark}{Remark}

\renewcommand{\baselinestretch}{0.99}
\begin{document}

\maketitle
\thispagestyle{empty}
\pagestyle{empty}

%%%%%%%%%%%%%%%%%%%%%%%%%%%%%%%%%%%%%%%%%%%%%%%%%%%%%%%%%%%%%%%%%%%%%%%%%%%%%%%%
\begin{abstract}
% A common belief is that distributed algorithms enable private computations of Optimal Power Flow (OPF) problem in power systems. They negate the need of sharing sensitive information, such as system loads, that is contained in algorithm sub-problems. This information, however, can still be inferred by a potential adversary from the responses of sub-problems to input signals across iterations. To provide formal privacy guarantees for distributed OPF computations, this paper leverages the framework of differential privacy and develops privacy-preserving distributed algorithms. The paper considers a distributed implementation of the OPF problem based on consensus Alternating Direction Method of Multipliers (ADMM) and proposes two methods to ensure differential privacy: dynamic and static perturbations of sub-problem solutions. The perturbations are random and drawn from a carefully parameterized Laplace distribution, so that the sensitive information cannot be reliably inferred by an adversary. The paper shows theoretically and numerically that the random perturbations provide a strong privacy protection for ADMM sub-problems with guarantees. Using a standard IEEE 118-node test case, the paper explores the fundamental trade-offs among privacy, algorithmic convergence, and optimality loss. 

Distributed algorithms enable private Optimal Power Flow (OPF)
computations by avoiding the need in sharing sensitive information
localized in algorithms sub-problems. However, adversaries can still
infer this information from the coordination signals exchanged across
iterations. This paper seeks formal privacy guarantees for distributed
OPF computations and provides differentially private algorithms for
OPF computations based on the consensus Alternating Direction Method
of Multipliers (ADMM). The proposed algorithms attain differential
privacy by introducing static and dynamic random perturbations of OPF
sub-problem solutions at each iteration. These perturbations are
Laplacian and designed to prevent the inference of sensitive
information, as well as to provide theoretical privacy guarantees for
ADMM sub-problems. Using a standard IEEE 118-node test case, the paper
explores the fundamental trade-offs among privacy, algorithmic
convergence, and optimality losses.

\end{abstract}
%%%%%%%%%%%%%%%%%%%%%%%%%%%%%%%%%%%%%%%%%%%%%%%%%%%%%%%%%%%%%%%%%%%%%%%%%%%%%%%%
\begingroup
\allowdisplaybreaks

\section{Introduction}

Centralized OPF computations operate over large datasets of system
parameters, such as electrical loads, and their unintended release
poses privacy risks for data owners. Recognizing these risks, the
literature suggests replacing the centralized computations with
distributed algorithms \cite{molzahn2017survey}, e.g., using the
well-known Alternating Direction Method of Multipliers (ADMM). These
algorithms distribute OPF computations among sub-problems that
coordinate through primal and dual coordination signals without
sharing the parameters used in their local computations and thus,
preserving privacy. However, in the presence of side information,
adversaries can reverse-engineer local parameters from observed
coordination signals \cite{shokri2017membership}. To overcome this
limitation, this paper augments these ADMM-based OPF algorithms with
{\it differential privacy}.

Differential privacy, first formalized by Dwork {\it et al.}
\cite{dwork2014algorithmic,dwork2006calibrating}, is a theoretical
framework {\it quantifying} the privacy risk associated with computing
functions (queries) on datasets with sensitive information. It ensures
that the same query applied to two adjacent datasets, i.e., differing
by one item, return essentially similar results (i.e., up to specified
parameters), thus preventing adversaries from learning any substantial
information over individual items. Chatzikokolakis {\it et al.}
\cite{chatzikokolakis2013broadening} generalized this concept to a
metric-based differential privacy for cases where publicly known participants
have sensitive data to protect. For instance, in power systems,
instead of hiding the presence of industrial customers, their
electrical loads may need to be obfuscated (up to a certain threshold)
to avoid revealing their commercial activities. Traditionally, differential privacy is
achieved by adding Laplacian noise to query outputs and the noise can
be calibrated using a small set of parameters (e.g., the privacy
loss $\varepsilon$) to control the differences between the outputs
on two adjacent datasets and obtain the guarantee known as
$\varepsilon-$differential privacy \cite{dwork2014algorithmic}.

% \noindent
{\it Contributions:} This paper applies differential privacy to
distributed OPF computations using a consensus ADMM, where the OPF
sub-problems coordinate voltage variables iteratively without
disclosing their local load parameters.  The paper first introduces an
adversarial model and shows that, in the presence of additional
information, adversaries can infer local load parameters from the
sub-problem responses to the coordination signals. To remedy this
privacy leak, the paper introduces two privacy-preserving ADMM
algorithms for OPF computations using static
(SP-ADMM) and  dynamic (DP-ADMM) random perturbations of sub-problem solutions across
iterations. The paper proves that the two algorithms ensure
$\varepsilon-$differential privacy but differ in the amount of noise
they introduce. The DP-ADMM algorithm ensures privacy at each
iteration, but needs to scale the noise magnitude in order to minimize
the privacy loss across multiple iterations. On the other hand, the
SP-ADMM preserves differential privacy across all iterations uniformly
but the worst-case sensitivity of its sub-problem solutions to load
datasets must be defined ahead of the algorithm iterations. Numerical
experiments highlight that, with a fine calibration of privacy
parameters, the inference of loads from primal-dual coordination
signals is equivalent to random guessing, {\it even if an adversary
  acquires all but one unknown sub-problem parameters}. The
experiments also explore the convergence properties of the two
algorithms and evaluate their fidelity with respect to the non-private
ADMM. In particular, despite similar privacy properties, DP-ADMM
results in smaller optimality losses, while SP-ADMM, exhibits faster
convergence rates.

% \vspace{1cm}
% Motivated by this scenario, the paper contributes by introducing privacy-preserving ADMM algorithms based on dynamic and static random perturbations of sub-problem solutions, termed DP-ADMM and SP-ADMM, respectively. The paper derives theoretical guarantees that the two algorithms provide $\varepsilon-$differential privacy, but the specification of the noise is different among the algorithms. The DP-ADMM preserves privacy for a single iteration of the algorithm. To bound the privacy loss across a finite number of iterations, the DP-algorithm extends to include the sequential composition result from \cite{dwork2014algorithmic}. The SP-ADMM preserves differential privacy across all iterations by design, but the worst-case sensitivity of sub-problem solutions to load datasets needs to be defined in advance of the algorithm iterations. A series of numerical experiments  illustrates that with a fine tuning of privacy parameters, the inference of loads from primal-dual coordination signals conforms to random guessing {\it even if an adversary acquires all but one unknown sub-problem parameters}. The experiments also analyze the convergence properties of the two algorithms and evaluate their fidelity with respect to the non-private results. In particular, despite similar privacy properties, it is shown that dynamic perturbations result in smaller optimality loss and power mismatch at the planning stage, both induced due to the noise in computations. The static perturbations, however, demonstrate stronger converge.

% \noindent
{\it Related work:} Differentially private distributed computation
was first introduced by Zhang {\it et al.}  \cite{zhang2016dynamic}
for the unconstrained empirical risk minimization (ERM) problem. The
authors distribute the ERM problem using ADMM and obtain differential
privacy for local training datasets by adding noise to either primal
or dual variables. The privacy guarantees, however, are provided for a
single ADMM iteration. Moreover, the results hold for the
unconstrained ERM problem, and are not appropriate for heavily
constrained OPF computations. Han {\it et al.}
\cite{han2016differentially} build a private distributed projected
gradient descent algorithm with gradient perturbations for electrical
vehicle charging, preventing the inference of charging power from
coordination signals. In the OPF context, Mak {\it et al.}
\cite{mak2019privacy} extend the centralized private release of OPF
test cases \cite{fioretto2018constrained,fioretto2019differential} to
a distributed ADMM-based algorithm. However, the work is meant for the
private release of input datasets and does not provide the OPF
solution itself.

\section{Towards Distributed OPF Computations}
Consider a power system as an undirected graph
$\Gamma\braceit{\set{B},\Lambda}$, where $\set{B}$ is the set of nodes
and $\Lambda$ is the set of transmission lines. Each transmission line
has a susceptance $\beta\in\mathbb{R}_{+}^{|\Lambda|}$ and a capacity
$\overline{f}\in\mathbb{R}_{+}^{|\Lambda|}$. The mapping functions
$\text{s}:\Lambda\mapsto\set{B}$ and $\text{r}:\Lambda\mapsto\set{B}$
are used to return the sending and receiving ends of lines,
respectively. The network topology is described by a weighted
Laplacian matrix $B\in\mathbb{R}^{|\set{B}|\times|\set{B}|}$, where the
weights on the lines are given by their susceptances. The network
loads are given by vector $d\in\mathbb{R}_{+}^{|\set{B}|}$. Each node
generates an amount $p\in\mathbb{R}^{|\set{B}|}$ of real power in the
interval $[\underline{p},\overline{p}]$ for a cost given by a
quadratic function whose second- and first-order coefficients are
$c_{2}\in\mathbb{R}_{+}^{|\set{B}|}$ and
$c_{1}\in\mathbb{R}_{+}^{|\set{B}|}$, respectively.  The OPF solution
amounts to generator set-points $p\in\mathbb{R}^{|\set{B}|}$ and
voltage angles $\theta\in\mathbb{R}^{|\set{B}|}$, obtained from
the optimization
\begin{subequations}\label{OPF_center}
\begin{align}
    \minimize{p,\theta} \quad& c\braceit{p}=\sum_{i\in\set{B}}c_{2i}p_{i}^2 + c_{1i}p_{i} \label{OPF_center_obj}\\
    \st \quad&\underline{p}_{i}\leqslant p_{i}\leqslant\overline{p}_{i},\;,\forall i\in\set{B},\label{OPF_center_gen_limit}\\
    &-\overline{f}_{l}\leqslant\beta_{\ell}\braceit{\theta_{\text{s}\braceit{\ell}}-\theta_{\text{r}\braceit{\ell}}}\leqslant\overline{f}_{l},\;\forall \ell\in\Lambda,\label{OPF_center_flow_limit}\\ 
    &\sum_{j\in\set{B}}B_{ij}\theta_{j} = p_{i} - d_{i},\forall i\in\set{B}, \label{OPF_center_balance}
\end{align}
\end{subequations}
which minimizes the total generation cost
\eqref{OPF_center_obj}. Inequality constraints
\eqref{OPF_center_gen_limit} and \eqref{OPF_center_flow_limit}
respectively ensure that generation and power flows are within their
corresponding limits. Equations \eqref{OPF_center_balance} ensure
the balance among the load, generation and power flow injection at
every network node.

The centralized computation in \eqref{OPF_center} requires that all
network parameters are submitted to a central entity. To avoid the
need for sharing network parameters, one may consider instead, a
distributed OPF computation, where the network is arbitrarily split
into zones $z\in\set{Z}$ \cite{erseghe2015distributed}. The domestic
nodes of each zone $z$ are collected in a set $\set{R}_{z}$, such that
$\set{R}_{z}\cap\set{R}_{z'}=\emptyset,\forall z \neq z'$. The
extended set $\set{V}_{z}$ contains the domestic nodes of, and
adjacent nodes to, zone $z$. The set of end nodes from the
transmission lines adjacent to zone $z$ is then defined as
$\set{M}_{z}=\set{V}_{z}\cap\set{V}_{z'}.\forall z \neq z'$. To enable
the distributed computation, the voltage angles are duplicated per
zone, i.e., they are redefined as
$\theta\in\mathbb{R}^{|\set{B}|\times|\set{Z}|}$. Towards the purpose, the
following consensus constraint is enforced:
\begin{align}
    \theta_{iz}=\overline{\theta}_{i}:\mu_{iz},\; \forall z\in\set{Z},\;\forall i\in\set{M}_{z},\label{consensus_constraint}
\end{align}
where $\theta_{iz}$ is a local copy of the voltage angle at node $i$,
$\overline{\theta}\in\mathbb{R}^{|\set{B}|}$ is the consensus variable,
and $\mu\in\mathbb{R}^{|\set{B}|\times|\set{Z}|}$ is the dual variable of
the consensus constraint. This decomposition separates the feasibility
region \eqref{OPF_center_gen_limit}-\eqref{OPF_center_balance} per
zone, so we denote the constraint set of each zone by
$\set{F}_{z}$. Now, the computation boils down to the optimization of
the following partial Lagrangian function:
\begin{align*}
    \maximize{\mu}\;\minimize{p,\theta,\overline{\theta}}\;&  \mathcal{L}(\mu,p,\theta,\overline{\theta})=c\braceit{p} + \mu_{z}^{\top}\braceit{\overline{\theta}-\theta_{z}}\\
    \st\quad&p,\theta\in\cap_{z\in\set{Z}}\set{F}_{z},
\end{align*}
where the objective function includes the dualized consensus constraint \eqref{consensus_constraint} with $\mu_{z}$ and $\theta_{z}$ being $z^{\text{th}}$ columns of $\mu$ and $\theta$, respectively. The distributed OPF computation is thus enabled by the following ADMM algorithm: 
\begingroup
\allowdisplaybreaks
\begin{align*}
    &\theta_{z}^{\kplus}\leftarrow\argmin{(p,\theta_{z})\in\set{F}_{z}}\;\mathcal{L}(\mu^{k},p,\theta_{z},\overline{\theta}^{k}) + \frac{\rho}{2}\norm{\overline{\theta}^{k}-\theta_{z}}_{2}^{2},\;\forall z\in\set{Z},\\
    &\overline{\theta}^{\kplus}\leftarrow\argmin{\overline{\theta}}\;\mathcal{L}(\mu^{k},\theta^{\kplus},\overline{\theta}) + \frac{\rho}{2}\sum_{z\in\set{Z}}\norm{\overline{\theta}-\theta_{z}^{\kplus}}_{2}^{2},\\
    &\mu_{z}^{\kplus}\leftarrow\mu_{z}^{k}+\rho\braceit{\overline{\theta}^{\kplus}-\theta_{z}^{\kplus}},\;\forall z\in\set{Z},
\end{align*}
\endgroup where $k$ is an iteration index and the squared norms denote
the ADMM regularization terms augmented with a non-negative penalty
factor $\rho.$ The algorithm is indeed distributed, as it solely
requires the exchange and update of primal and dual coordination signals between the neighboring zones. 

The rest of the paper makes the following
assumption.
\begin{assumption}
\label{base_assumption}
The function $c(p)$ is convex and strictly monotone in $p$, the set $\set{F}_{z}$ is compact and convex for all $z\in\set{Z}$,  and $\cap_{z\in\set{Z}}\set{F}_{z}$ has a non-empty interior. As a result, the distributed OPF algorithm converges to a unique optimal solution in a finite number of iterations \cite{boyd2011distributed}.  
\end{assumption}

\section{Differential Privacy for OPF}

The privacy goal in this paper is to ensure that individual loads cannot be inferred from the outputs of the ADMM sub-problems. Each sub-problem is thus considered as a {\it query} $\mathcal{Q}_{z}^{k}$ that maps the dataset of loads $\set{D}_{z}=\{d_{i}\}_{i\in\set{R}_{z}}$ to a vector of all voltage angles $\{\theta_{iz}^{\kplus}\}_{i\in\set{M}_{z}}$ to be released at iteration $k$. Under Assumption \ref{base_assumption}, each sub-problem has a unique response for a given load dataset, i.e., it computes a one-to-one mapping of the load dataset to the voltage solution. Hence, the release of the voltage solution leads to the leakage of the load dataset. 

The dependencies between sub-problem solutions and their load datasets can be weakened by making queries $\mathcal{Q}_{z}^{k}$ differentially private, i.e., by adding a carefully calibrated noise to their outputs. More precisely, the added noise aims at making
the outputs for two {\em adjacent} load datasets $\set{D}_{z}$ and $\set{D}_{z}^{\prime}$ indistinguishable from each other. 

\begin{definition}[Adjacency \cite{chatzikokolakis2013broadening}]
$\set{D}_{z}=\{d_{i}\}_{i\in\set{R}_{z}}$ and $\set{D}_{z}^{\prime}=\{d_{i}^{\prime}\}_{i\in\set{R}_{z}}$  are $\alpha-$adjacent datasets,
denoted by $\set{D}_{z}\sim_{\alpha}\set{D}_{z}^{\prime}$, if they differ in one element by $\alpha$, i.e.,
\begin{align*}
\exists i \textrm{~s.t.~} \norm{d_i - {d}_{i}^{\prime}}_{1}\leqslant \alpha \;\land\; d_j = {d}_{j}^{\prime}, \forall j \neq i.
\end{align*}
\end{definition}

Differential privacy relies on the concept of global sensitivity to calibrate the noise. 
\begin{definition}[Global Query Sensitivity]\label{def_global_sensitivity}
The global sensitivity $\Delta_{z}$ of query $\mathcal{Q}^k_{z}$ is defined by
\begin{align*}
    \Delta_{z}^{k}\coloneqq\maximize{\set{D}_{z}\sim_{\alpha}\set{D}_{z}^{\prime}}\;\norm{
    \mathcal{Q}_{z}^k(\set{D}_{z})-\mathcal{Q}_{z}^{k}(\set{D}_{z}^{\prime})
    }_{1}.
\end{align*}
where $\set{D}_{z}$ and $\set{D}_{z}^{\prime}$ belong to the universe
$\mathcal{D}_z$ of all datasets of interest for zone $z$. Note that
the datasets in $\mathcal{D}_{z}$ are projections of the globally
feasible solutions of interest.
\end{definition}

%We consider that $\mathcal{D}_{z}$ includes all datasets that return feasible solutions of sub-problems and that %$\cap_{z\in\set{Z}}\set{F}_{z}\braceit{\set{D}_{z}}$ for all $\set{D}_{z}\in\mathcal{D}_{z}$ has a non-empty interior.\noindent

In practice, especially in off-peak hours, the global sensitivity is
overly pessimistic. The notion of local query sensitivity can be used
to obtain more precise upper bounds.
\begin{definition}[Local Query Sensitivity]\label{def_local_sensitivity}
  The local query sensitivity of query $\mathcal{Q}_{z}^{k}$ with
  respect to $\set{D}_{z}$, denoted by
  $\delta_{z}(\set{D}_{z})$, is defined as
\begin{align}\label{local_sens}
    \delta_{z}(\set{D}_{z}) = \maximize{\set{D}^{\prime} \sim_{\alpha} \set{D}_{z}}\;\norm{
    \mathcal{Q}_{z}^{k}(\set{D}_{z})-\mathcal{Q}_{z}^{k}(\set{D}_{z}^{\prime})
    }_{1}.
\end{align}
\end{definition}
\noindent
For simplicity, when $\set{D}_{z}$ is clear from the context, $\delta_{z}$ is used to denote $\delta_{z}(\set{D}_{z})$.
\begin{remark}
The maximal local sensitivity of $\mathcal{Q}_{z}^{k}$ depends not
only on the magnitude of loads $\{d_{i}\}_{i\in\set{R}_{z}}$, but also
on their relative position with respect to the nodes in $\set{M}_{z}$.
\end{remark}

This work introduces noise drawn from a zero-mean Laplace distribution
with scale $b$, denoted by $\text{Lap}(b)$ for short, with a
probability density function
$\text{Lap}(\xi|b)=\frac{1}{2b}\text{exp}(-\frac{\norm{\xi}_{1}}{b})$.
It can be used to attain differential privacy for numerical queries as
per the following result \cite{dwork2014algorithmic}.
\begin{theorem}[Laplace mechanism]\normalfont
Let $\mathcal{Q}:\set{D}\mapsto\mathbb{R}$ be a query that maps dataset $\set{D}$ to real numbers, and let $\Delta$ be a query sensitivity. The Laplace mechanism that outputs $\mathcal{Q}\braceit{\set{D}}+\xi$ with $\xi\sim\text{Lap}\braceit{\Delta/\varepsilon}$ is $\varepsilon-$differential privacy, i.e.,
\begin{align*}
    \mathbb{P}[\mathcal{Q}(\set{D})+\xi]\leqslant\mathbb{P}[\mathcal{Q}(\set{D}^{\prime})+\xi]\text{exp}(\varepsilon),
\end{align*}
where $\set{D}\sim_{\alpha}\set{D}^{\prime}$ are any $\alpha-$adjacent datasets. 
\end{theorem}

The parameter $\varepsilon$, called the {\it privacy loss}, bounds the
multiplicative difference between distributions of query outputs on
any two $\alpha$-adjacent datasets. Stronger privacy requirements can
be obtained by choosing smaller values for $\varepsilon$ and larger
values for $\alpha$. The last building block is the {\it sequential
  composition} theorem \cite{dwork2014algorithmic}, which
characterizes the guarantees for sequential applications of
differential privacy.
\begin{theorem}[Sequential composition]\label{th:compositon}
Consider $T$ runs of query function $\{\mathcal{Q}_{z}^{t}(\set{D}_{z})\}_{t=1}^{T}$ such that every run depends on the result of the previous runs, i.e., 
\begin{align*}
    \mathcal{Q}_{z}^{t}\braceit{\set{D}_{z}} =  \mathcal{Q}_{z}^{t}\braceit{\set{D}_{z},\mathcal{Q}_{z}^{1}(\set{D}_{z}),\mathcal{Q}_{z}^{2}(\set{D}_{z}),\dots,\mathcal{Q}_{z}^{t-1}(\set{D}_{z})}.
\end{align*}
Suppose that $\mathcal{Q}_{z}^{t}$ preserve $\varepsilon_{t}-$differential privacy for $\mathcal{Q}_{z}^{t'}$ for all $t'<t$. Then, the $T$-tuple mechanism $\mathcal{Q}_{z}=\braceit{\mathcal{Q}_{z}^{1},\mathcal{Q}_{z}^{2},\dots,\mathcal{Q}_{z}^{t}}$ preserves $\sum_{t=1}^{T}\varepsilon_{t}-$differential privacy.  
\end{theorem}

\section{Private Distributed OPF Algorithms}

\begin{algorithm}[!t]
{
\setlength{\abovedisplayskip}{3pt}
\setlength{\belowdisplayskip}{3pt}
\small{
  \caption{The SP-ADMM Algorithm}
  \label{alg:SDP-ADMM}
\begin{algorithmic}[1]
% \Procedure{Euclid}{$a,b$}\Comment{The g.c.d. of a and b}
\State {\bf Input:} Datasets $\set{D}_{z}$, privacy parameters $\varepsilon, \alpha,$ algorithmic parameters $\gamma,\rho,K,\overline{\theta}^{1},\mu^{1}$
\State Draw random samples $\xi_{z}\sim\text{Lap}(\frac{\Delta_{z}}{\varepsilon}),\forall z\in\set{Z}$
\While{$k\not=K$ or $\sum_{z\in\set{Z}}\norm{\tilde{\theta}_{z}^{k}-\overline{\theta}^{\kplus}}_{2}\leqslant\gamma$}
\State Update voltage angles $\theta_{z}^{\kplus},\forall z\in\set{Z},$ by solving $$\minimize{(p,\theta_{z})\in\set{F}_{z}}\mathcal{L}_{z}(\mu^{k},p,\theta_{z},\overline{\theta}^{k}) + \frac{\rho}{2}\norm{\overline{\theta}^{k}-\theta_{z}}_{2}^{2}$$
\State Perturb sub-problem solutions $\tilde{\theta}_{z}^{\kplus}=\theta_{z}^{\kplus}+\xi_{z},\forall z\in\set{Z},$
\State Update consensus variables $\overline{\theta}_{i}^{\kplus},\forall i\in\set{M}_{z},z\in\set{Z},$ as
$$\minimize{\overline{\theta}}\;\mathcal{L}(\mu^{k},\tilde{\theta}^{\kplus},\overline{\theta}) + \frac{\rho}{2}\sum_{z\in\set{Z}}\norm{\overline{\theta}-\tilde{\theta}_{z}^{\kplus}}_{2}^{2}$$
\State Update dual variables $\mu_{z}^{\kplus},\forall z\in\set{Z},$ by solving 
$$\mu_{z}^{\kplus}\leftarrow\;\mu_{z}^{k}+\rho\braceit{\overline{\theta}_{z}^{\kplus}-\tilde{\theta}_{z}^{\kplus}}$$
\State Iteration update $k\leftarrow k + 1$
\EndWhile\label{euclidendwhile}
\State \textbf{Output:} Private OPF solution.
% \EndProcedure
\end{algorithmic}
}
}
\end{algorithm}

This section presents two differentially private ADMM algorithms for distributed OPF computations. 

\subsection{The SP-ADMM Algorithm}

The SP-ADMM algorithm relies on a key insight: {\em the upper bound on the global query sensitivity is independent from the input coordination signals}.

\begin{proposition}\label{prop_global_sensitivity}
The global sensitivity of the ADMM sub-problem $\mathcal{Q}_{z}^{k}$ is upper-bounded by $\text{max}_{\set{D}_{z}\in\mathcal{D}_{z}}\text{max}_{i}\{d_{i}\}_{i\in\set{R}_{z}}$.
\end{proposition}
\begin{proof}
The result follows from the nodal balance constraint \eqref{OPF_center_balance}. Consider that the voltage angle sensitivity to load changes reduces with the size of the network graph. Therefore, the worst-case sensitivity is observed in a minimal size two-node network, i.e.,  
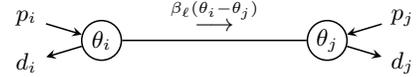
\begin{center}
\tikzstyle{NODE} = [circle, minimum width=15,text centered, draw=black, fill=white,line width=0.2mm]
\begin{tikzpicture}[font=\small]
\node [align=center] (node_i) [NODE,inner sep=0.5pt] {$\theta_{i}$};
\node [align=center] (node_j) [NODE,inner sep=0.5pt,right of = node_i, xshift = 2cm] {$\theta_{j}$};
\draw [black,line width = 0.02cm] (node_i) -- node[midway,above]{$\overset{\beta_{\ell}(\theta_{i}-\theta_{j})}{\longrightarrow}$} (node_j);
\node [align=center] (p_i) [left of = node_i, yshift=0.3cm] {$p_{i}$};
\node [align=center] (d_i) [left of = node_i, yshift=-0.3cm] {$d_{i}$};
\draw [black,line width = 0.02cm,->,>=stealth] (p_i) --  (node_i);
\draw [black,line width = 0.02cm,<-,>=stealth] (d_i) --  (node_i);
\node [align=center] (p_j) [right of = node_j, yshift=0.3cm] {$p_{j}$};
\node [align=center] (d_j) [right of = node_j, yshift=-0.3cm] {$d_{j}$};
\draw [black,line width = 0.02cm,->,>=stealth] (p_j) --  (node_j);
\draw [black,line width = 0.02cm,<-,>=stealth] (d_j) --  (node_j);
\end{tikzpicture}
\end{center}
where node $i$ is chosen as a reference node, i.e., $\theta_{i}=0.$ Then, the power balance at node $j$ is 
$-\beta_{\ell}\theta_{j}=p_{j}-d_{j}$. The worst-case sensitivity is provided given that $p_{j}=0$, so the change of load directly translates into $\beta_{\ell}\theta_{j}$. As $\beta_{\ell}\gg1$ in power system networks, the change of $\theta_{j}$ is upper-bounded by the magnitude of $d_{j}$. In turn, $d_{j}$ has to be chosen as the largest feasible load in a dataset universe $\mathcal{D}_{z},\forall z\in\set{Z}$.   
%   The result follows directly from Equation \eqref{OPF_center_balance} and the fact that $\beta_{\ell}\gg1,\forall \ell\in\Lambda.$
\end{proof}

As a result, the SP-ADMM algorithm, shown in Algorithm
\ref{alg:SDP-ADMM}, generates the Laplacian noise $\xi_{z}\in\mathbb{R}^{|\set{M}_{z}|},\forall z\in\set{Z},$  once, at the
beginning of the algorithm, using an upper bound $\Delta_z$ on
the global sensitivity. It takes the dataset, privacy and algorithm parameters  as inputs, and runs ADMM
iterations until reaching iteration limit $K$ or the primal residual
is below the tolerance $\gamma$. Unlike the conventional ADMM, once a
sub-problem produces the optimal response to the dual and consensus
variables, the algorithm perturbs the response with the initially generated noise. The perturbed solution $\tilde{\theta}_{z}^{\kplus}$
then participates in the consensus and dual variable updates.

\subsection{The DP-ADMM Algorithm}
\begin{algorithm}[t]
{
\setlength{\abovedisplayskip}{3pt}
\setlength{\belowdisplayskip}{3pt}
\small{
\caption{The DP-ADMM Algorithm}\label{alg:DDP-ADMM}
\begin{algorithmic}[1]
% \Procedure{Euclid}{$a,b$}\Comment{The g.c.d. of a and b}
\State {\bf Input:} Datasets $\set{D}_{z}$, privacy parameters $\varepsilon, \alpha,$ algorithmic parameters $\gamma,\rho,K,\overline{\theta}^{1},\mu^{1}$
\While{$k\not=K$ or $\sum_{z\in\set{Z}}\norm{\tilde{\theta}_{z}^{k}-\overline{\theta}^{\kplus}}_{2}\leqslant\gamma$}
\State Update voltage angles $\theta_{z}^{\kplus},\forall z\in\set{Z},$ by solving $$\minimize{(p,\theta_{z})\in\set{F}_{z}}\mathcal{L}_{z}(\mu^{k},p,\theta,\overline{\theta}^{k}) + \frac{\rho}{2}\norm{\overline{\theta}^{k}-\theta_{z}}_{2}^{2}$$
\State For $\mu_{z}^{k}$ and $\overline{\theta}^{k}$, compute sensitivity $\delta_{z}^{k},\forall z\in\set{Z},$ by solving
\begin{align*}
    \delta_{z}^{k}=\maximize{\set{D}^{\prime}\in\mathcal{D}}\;&\norm{
    \mathcal{Q}_{z}^{k}(\set{D}_{z})-\mathcal{Q}_{z}^{k}(\set{D}_{z}^{\prime})
    }_{1},\\
    \st\;&\norm{\set{D}_{z}-\set{D}_{z}^{\prime}}_{1}\leqslant\alpha
\end{align*}
\State Perturb sub-problem solutions by $\xi_{z}^{k}\sim\text{Lap}(\frac{\delta_{z}^{k}}{\varepsilon}),\forall z\in\set{Z},$ $$\tilde{\theta}_{z}^{\kplus}=\theta_{z}^{\kplus}+\xi_{z}^{k}$$
\State Update consensus variables $\overline{\theta}_{i}^{\kplus},\forall i\in\set{M}_{z},z\in\set{Z},$ as
$$\minimize{\overline{\theta}}\;\mathcal{L}(\mu^{k},\tilde{\theta}^{\kplus},\overline{\theta}) + \frac{\rho}{2}\sum_{z\in\set{Z}}\norm{\overline{\theta}_z-\tilde{\theta}_{z}^{\kplus}}_{2}^{2}$$
\State Update dual variables $\mu_{z}^{\kplus},\forall z\in\set{Z},$ by solving 
$$\mu_{z}^{\kplus}\leftarrow\;\mu_{z}^{k}+\rho\braceit{\overline{\theta}_{z}^{\kplus}-\tilde{\theta}_{z}^{\kplus}}$$
\State Iteration update $k\leftarrow k+1$
\EndWhile\label{euclidendwhile}
\State \textbf{Output:} Private OPF solution.
% \EndProcedure
\end{algorithmic}
}
}
\end{algorithm}

The DP-ADMM algorithm is fundamentally different: it uses the concept of
local query sensitivity to perturb the phase angles differently at each
iteration. Its key insight is the recognition that {\em the local
sensitivity of the queries/sub-problems can be obtained by solving an
optimization problem.} As a result, for each iteration, the DP-ADMM
perturbs the phase angles using the local sensitivity. The DP-ADMM is
outlined in Algorithm \ref{alg:DDP-ADMM}.  The noise $\xi_{z}^{k}\in\mathbb{R}^{|\set{M}_{z}|},\forall z\in\set{Z},$ is
dynamically updated respecting the change of local sensitivity
$\delta_{z}^{k}$ on $\alpha-$adjacent datasets. The sensitivity is
obtained by identifying the individual load, whose $\alpha-$change
brings the maximal change of the sub-problem solution. The rest of the
algorithm is similar to SP-ADMM.

\subsection{Properties}

This section reviews the properties of the two algorithms. 

\begin{theorem}\label{th:DP_privacy}
{
\setlength{\abovedisplayskip}{1.5pt}
\setlength{\belowdisplayskip}{1.5pt}
    Let $\tilde{\mathcal{Q}}_{z}^{k}(\set{D}_{z})$ be a randomized sub-problem of zone $z$ acting on optimization dataset $\set{D}_{z}$ i.e., $$\tilde{\mathcal{Q}}_{z}^{k}(\set{D}_{z})=\theta_{z}^{k+1}+\xi_{z}^{k}.$$
    Let $\delta_{z}^{k}$ be a sensitivity of $\mathcal{Q}_{z}^{k}(\set{D}_z)$ for all $\alpha-$adjacent dataset $D_{z}^{\prime}$. Then, if the random perturbation $\xi_{z}^{k}$ is sampled from the probability distribution with density function $\text{Lap}(\frac{\delta_{z}^{k}}{\varepsilon})$, then $\tilde{\mathcal{Q}}_{z}^{k}(\set{D})$ provides  $\varepsilon-$differential privacy at iteration $k$, i.e., 
    $$\mathbb{P}[\tilde{\mathcal{Q}}_{z}^{k}(\set{D}_{z})]\leqslant\mathbb{P}[\tilde{\mathcal{Q}}_{z}^{k}(\set{D}_{z}^{\prime})]\text{exp}(\varepsilon),\;\forall \set{D}_{z}^{\prime}\in\mathcal{D}_{z}.$$
}
\end{theorem}
\begin{proof}
The proof follows a similar line of arguments as for numerical queries. We consider the ratio of probabilities that the query $\tilde{\mathcal{O}}_{z}^{k}$ returns the same solution $\hat{\theta}_{z}^{k+1}$ on two $\alpha-$adjacent datasets $\set{D}_{z}\sim_{\alpha}\set{D}_{z}'$ at ADMM iteration $k$:
\begin{align}\label{quick_ref_1}
    &\frac{\mathbb{P}[\tilde{\mathcal{Q}}_{z}^{k}(\set{D}_{z})\in\hat{\theta}_{z}^{k+1}]}{\mathbb{P}[\tilde{\mathcal{Q}}_{z}^{k}(\set{D}_{z}^{\prime})\in\hat{\theta}_{z}^{k+1}]}
    =\frac{\mathbb{P}[\mathcal{Q}_{z}^{k}(\set{D}_{z})+\text{Lap}(\xi_{z}|\frac{\delta_{z}^{k}}{\varepsilon})\in\hat{\theta}_{z}^{k+1}]}{\mathbb{P}[\mathcal{Q}_{z}^{k}(\set{D}_{z}^{\prime})+\text{Lap}(\xi_{z}|\frac{\delta_{z}^{k}}{\varepsilon})\in\hat{\theta}_{z}^{k+1}]}\nonumber\\
    &=
    \prod_{{i\in\set{M}_{z}}}\frac{\frac{\varepsilon}{2\delta_{z}^{k}}\text{exp}\left(-\frac{\norm{\xi_{iz}-\mathcal{Q}_{iz}^{k}(\set{D}_{z})}_{1}}{\delta_{z}^{k}}\right)}{\frac{\varepsilon}{2\delta_{z}^{k}}\text{exp}\left(-\frac{\norm{\xi_{iz}-\mathcal{Q}_{iz}^{k}(\set{D}_{z}^{\prime})}_{1}}{\delta_{z}^{k}}\right)}\nonumber\\
    &= 
     \prod_{{i\in\set{M}_{z}}}\text{exp}\braceit{\frac{\varepsilon\braceit{\norm{\xi_{iz}-\mathcal{Q}_{iz}^{k}(\set{D}_{z}^{\prime})}_{1} - \norm{\xi_{iz}-\mathcal{Q}_{iz}^{k}(\set{D}_{z})}_{1}}}{\delta_{z}^{k}}}\nonumber\\
    &\leqslant\prod_{{i\in\set{M}_{z}}}\text{exp}\braceit{\frac{\varepsilon\norm{\mathcal{Q}_{iz}^{k}(\set{D}_{z})-\mathcal{Q}_{iz}^{k}(\set{D}_{z}^{\prime})}_{1}}{\delta_{z}^{k}}}\nonumber\\
    &=\text{exp}\braceit{\frac{\varepsilon\norm{\mathcal{Q}_{z}^{k}(\set{D}_{z})-\mathcal{Q}_{z}^{k}(\set{D}_{z}^{\prime})}_{1}}{\delta_{z}^{k}}},
\end{align}
where the second equality follows from the definition of the probability
density function of the Laplace distribution, and the inequality
follows from the inequality of norms. Recall Definition
\ref{def_local_sensitivity} of the local sensitivity. Hence, by
substituting \eqref{local_sens} in \eqref{quick_ref_1}, we obtain the
desired result.
\end{proof}

\begin{remark}
Theorem \ref{th:DP_privacy} holds not only for the DP-ADMM, but also for the SP-ADMM algorithm when used with an upper bound on the global sensitivity.
\end{remark}

\begin{remark}
Theorem \ref{th:DP_privacy} makes use of the local sensitivity
$\delta_{z}^{k}$, thus attaining local $\varepsilon-$differential
privacy. By substituting $\delta_{z}^{k}$ with the global query
sensitivity $\Delta_{z}$, the algorithms provide global
$\varepsilon-$differential privacy. The robustness of the two
approaches to privacy attacks is analyzed in Section
\ref{experiments}.
\end{remark}

Observe that every new iteration of the DP-ADMM algorithm reveals more information to an adversary, thus diminishing the privacy guarantee. Assume that the algorithm implementation can
limit the adversary to observing $T$ iterations, e.g., by using
secure switching of communication channels. Then, the following
result applies. 

\begin{theorem}\label{th:compositon_two}
Let
$\tilde{\mathcal{Q}}_{z}^{k}(\set{D}_{z})=\theta_{z}^{k+1}+\xi_{z}^{k}$
be a randomized query as specified in Theorem \ref{th:DP_privacy} with
the difference that the noise $\xi_{z}^{k}$ is drawn from
$\text{Lap}(T\frac{\delta_{z}^{k}}{\varepsilon}).$ Then, DP-ADMM
preserves $\varepsilon-$differential privacy across $T$ iterations.
\end{theorem}
\proof{It follows from combining Theorems \ref{th:compositon} and \ref{th:DP_privacy}}. 

Finally, observe that the feasibility of the OPF solution is not
affected by either dynamic or static perturbations, as the two
algorithms add noise only to the unconstrained consensus and dual
variable updates.

\section{Adversarial Problem}

The strength of differentially private algorithms is their robustness
to {\it side} information. The framework guarantees that, even if an
adversary obtains information on {\it all  but one} items in a dataset,
the privacy of the remaining one item is ensured. This section presents an
adversarial problem for this worst-case scenario of privacy attack on
OPF sub-problems.

Consider a set $\set{T}=\{k-T,\dots,k\}$ of ADMM iterations
observed by an adversary. Let $\tilde{\theta}_{z}^{\tplus}$ be a
response of each sub-problem $z\in\set{Z}$ to dual and consensus
variables $\mu_{z}^{t}$ and $\overline{\theta}_{z}^{t}$ at iteration
$t\in\set{T}.$ For sub-problem $z$, the adversarial inference problem
can be formulated as the following empirical risk minimization problem
across $T$ iterations:
\begin{subequations}\label{OPF_attack_model}
\begin{align}
    &\minimize{\hat{p}^{t},\hat{\theta}_{z}^{t},\hat{d}_{i}}\;\sum_{t\in\set{T}}c_{z}\braceit{\hat{p}^{t}} - [\mu_{z}^{t}]^{\top}\hat{\theta}_{z}^{t} \nonumber\\
    &\quad\quad\quad+\sum_{t\in\set{T}}\frac{\rho}{2}\norm{\overline{\theta}^{t}-\hat{\theta}_{z}^{t}}_{2}^{2} \nonumber\\
    &\quad\quad\quad\quad\quad\quad\quad+ \Upsilon\sum_{t\in\set{T}}\norm{\hat{\theta}_{z}^{t}-\tilde{\theta}_{z}^{\tplus}}_{2}\\
    \st\;&\text{Equations}\;\eqref{OPF_center_gen_limit}-\eqref{OPF_center_flow_limit}, \;\forall t\in\set{T},\\
    % &B\hat{\theta}_{z}^{t} = \hat{p}^{t} - [d_{1},\dots,\hat{d}_{i},\dots, d_{|\set{R}_{z}|}]^{\top},\;\forall t\in\set{T},\\
    &\sum_{m\in\set{V}_{z}}B_{nm}\hat{\theta}_{mz}^{t}=\hat{p}_{n}^{t}-d_{n}, \forall n\in\set{R}_{z}\backslash i, t\in\set{T},\\
    &\sum_{m\in\set{V}_{z}}B_{im}\hat{\theta}_{mz}^{t}=\hat{p}_{i}^{t}-\hat{d}_{i}, \forall t\in\set{T},
\end{align}
\end{subequations}
where decision variables are indicated with a $(\hat{\cdot})$
notation, and the rest are the parameters available to an
adversary. The unknown load magnitude $\hat{d}_{i}$ at node $i$ of interest is modelled as a
decision variable. An adversary seeks the value of $\hat{d}_{i}$ that
minimizes the Euclidean distance between the voltage variables
$\hat{\theta}_{z}^{\tplus}$ modeled in the adversarial problem and the
voltage solution $\tilde{\theta}_{z}^{\tplus}$ released by the
sub-problem at all iterations $t\in\set{T}$. By penalizing the
distance with a sufficiently large coefficient $\Upsilon$, an
adversary identifies the load magnitude that produces the same voltage
solution as that released by the sub-problem, thus identifying the
unknown load magnitude.

%%%%%%%%%%%%%%%%%%%%%%%%%%%%%%%%%%%%%%%%%%%%%%%%%%%%%%%%
% \iffalse
%%%%%%%%%%%%%%%%%%%%%%%%%%%%%%%%%%%%%%%%%%%%%%%%%%%%%%%%

\section{Numerical Experiments}\label{experiments}

This section examines the proposed Algorithms \ref{alg:DDP-ADMM} and
\ref{alg:SDP-ADMM} using a standard IEEE 118-node test case with a
3-zone lay-out taken from \cite[case 118-3]{guo2015intelligent}. The
algorithms are compared in terms of their robustness to privacy
attacks, convergence properties, and fidelity with respect to the
non-private ADMM algorithm. By default, we set ADMM penalty factor
$\rho=100$, iteration limit $K=300$, algorithm tolerance $\gamma=0.5$,
and coefficient $\Upsilon=10^{6}$. The privacy requirements are
selected such that the privacy loss is fixed $\varepsilon=1$ whereas
the adjacency coefficient varies in the range
$\alpha=\{1,2.5,5,7,10\}\%$. For the given algorithmic parameters, the
standard non-private ADMM converges to the optimal OPF solution in 59
iterations.

\subsubsection{Robustness to the Privacy Attacks} 

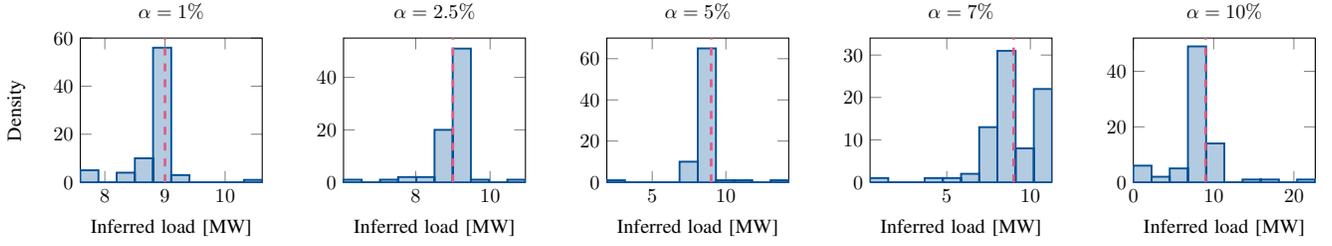
\begin{figure*}[t]
    \centering
    \begin{tikzpicture}[scale=1, every node/.style={scale=0.75}]
        \begin{axis}[
        xshift = 0cm,
        width=4cm,
        height=3.5cm,
        enlargelimits=false,
        xlabel = {Inferred load [MW]},
        ylabel = {Density},
        title = {$\alpha=1\%$},
        ylabel style = {yshift = -0.1cm},
        legend style={draw=none,fill=none,font=\scriptsize},
        legend cell align={left},
        ]
        \addplot[blue,fill=blue!30,hist,line width=0.025cm,opacity=1] table [y=a_1]{results_plot_data/hist_DDP.dat};
        \addplot[mark=none,line width=0.035cm,color=red,dashed] plot coordinates {(9,0) (9,60)};
        \end{axis}
        \begin{axis}[
        xshift = 3.5cm,
        width=4cm,
        height=3.5cm,
        enlargelimits=false,
        xlabel = {Inferred load [MW]},
        title = {$\alpha=2.5\%$},
        ylabel style = {yshift = -0.4cm},
        legend style={draw=none,fill=none,font=\scriptsize},
        legend cell align={left},
        ]
        \addplot[blue,fill=blue!30,hist,line width=0.025cm,opacity=1] table [y=a_2.5]{results_plot_data/hist_DDP.dat};
        \addplot[mark=none,line width=0.035cm,color=red,dashed] plot coordinates {(9,0) (9,55)};
        \end{axis}
        \begin{axis}[
        xshift = 7cm,
        width=4cm,
        height=3.5cm,
        enlargelimits=false,
        xlabel = {Inferred load [MW]},
        title = {$\alpha=5\%$},
        ylabel style = {yshift = -0.4cm},
        legend style={draw=none,fill=none,font=\scriptsize},
        legend cell align={left},
        ]
        \addplot[blue,fill=blue!30,hist,line width=0.025cm,opacity=1] table [y=a_5]{results_plot_data/hist_DDP.dat};
        \addplot[mark=none,line width=0.035cm,color=red,dashed] plot coordinates {(9,0) (9,70)};
        \end{axis}
        \begin{axis}[
        xshift = 10.5cm,
        width=4cm,
        height=3.5cm,
        enlargelimits=false,
        xlabel = {Inferred load [MW]},
        title = {$\alpha=7\%$},
        ylabel style = {yshift = -0.4cm},
        legend style={draw=none,fill=none,font=\scriptsize},
        legend cell align={left},
        ]
        \addplot[blue,fill=blue!30,hist,line width=0.025cm,opacity=1] table [y=a_7]{results_plot_data/hist_DDP.dat};
        \addplot[mark=none,line width=0.035cm,color=red,dashed] plot coordinates {(9,0) (9,34)};
        \end{axis}
        \begin{axis}[
        xshift = 14cm,
        width=4cm,
        height=3.5cm,
        enlargelimits=false,
        xlabel = {Inferred load [MW]},
        title = {$\alpha=10\%$},
        ylabel style = {yshift = -0.4cm},
        legend style={draw=none,fill=none,font=\scriptsize},
        legend cell align={left},
        ]
        \addplot[blue,fill=blue!30,hist,line width=0.025cm,opacity=1] table [y=a_10]{results_plot_data/hist_DDP.dat};
        \addplot[mark=none,line width=0.035cm,color=red,dashed] plot coordinates {(9,0) (9,52)};
        \end{axis}
        
    \end{tikzpicture}
    \vspace{-0.2cm}

    \caption{DP-ADMM: Results of privacy attack on the load sited at node 20. The plots show the densities of inferred load by an adversary across iterations for different adjacency coefficients for a single simulation run.}
    \label{fig:DP_ADMM_attack}
\end{figure*}

\begin{figure*}[t]
    \centering
    \begin{tikzpicture}[scale=1, every node/.style={scale=0.75}]
        \begin{axis}[
        xshift = 0cm,
        width=4cm,
        height=3.5cm,
        enlargelimits=false,
        xlabel = {Inferred load [MW]},
        ylabel = {Density},
        title = {$\alpha=1\%$},
        ylabel style = {yshift = -0.1cm},
        legend style={draw=none,fill=none,font=\scriptsize},
        legend cell align={left},
        ]
            \addplot[red,fill=red!30,hist,line width=0.025cm,opacity=1] table [y=x1]{results_plot_data/save_SDP_inference_results_global_sense.dat};
            \addplot[blue,fill=blue!30,hist,line width=0.025cm,opacity=0.75] table [y=x1]{results_plot_data/save_SDP_inference_results.dat};
            \addplot[mark=none,line width=0.035cm,color=red,dashed] plot coordinates {(9,0) (9,400)};
        \end{axis}
        \begin{axis}[
        xshift = 3.5cm,
        width=4cm,
        height=3.5cm,
        enlargelimits=false,
        xlabel = {Inferred load [MW]},
        title = {$\alpha=2.5\%$},
        ylabel style = {yshift = -0.4cm},
        legend style={draw=none,fill=none,font=\scriptsize},
        legend cell align={left},
        ]
            \addplot[red,fill=red!30,hist,line width=0.025cm,opacity=1] table [y=x2]{results_plot_data/save_SDP_inference_results_global_sense.dat};
            \addplot[blue,fill=blue!30,hist,line width=0.025cm,opacity=0.75] table [y=x2]{results_plot_data/save_SDP_inference_results.dat};
            \addplot[mark=none,line width=0.035cm,color=red,dashed] plot coordinates {(9,0) (9,400)};
        \end{axis}
        \begin{axis}[
        xshift = 7cm,
        width=4cm,
        height=3.5cm,
        enlargelimits=false,
        xlabel = {Inferred load [MW]},
        title = {$\alpha=5\%$},
        ylabel style = {yshift = -0.4cm},
        legend style={draw=none,fill=none,font=\scriptsize},
        legend cell align={left},
        ]
            \addplot[red,fill=red!30,hist,line width=0.025cm,opacity=1] table [y=x3]{results_plot_data/save_SDP_inference_results_global_sense.dat};
            \addplot[blue,fill=blue!30,hist,line width=0.025cm,opacity=0.75] table [y=x3]{results_plot_data/save_SDP_inference_results.dat};
            \addplot[mark=none,line width=0.035cm,color=red,dashed] plot coordinates {(9,0) (9,420)};
        \end{axis}
        \begin{axis}[
        xshift = 10.5cm,
        width=4cm,
        height=3.5cm,
        enlargelimits=false,
        xlabel = {Inferred load [MW]},
        title = {$\alpha=7\%$},
        ylabel style = {yshift = -0.4cm},
        legend style={draw=none,fill=none,font=\scriptsize},
        legend cell align={left},
        ]
            \addplot[red,fill=red!30,hist,line width=0.025cm,opacity=1] table [y=x4]{results_plot_data/save_SDP_inference_results_global_sense.dat};
            \addplot[blue,fill=blue!30,hist,line width=0.025cm,opacity=0.75] table [y=x4]{results_plot_data/save_SDP_inference_results.dat};
            \addplot[mark=none,line width=0.035cm,color=red,dashed] plot coordinates {(9,0) (9,400)};
        \end{axis}
        \begin{axis}[
        xshift = 14cm,
        width=4cm,
        height=3.5cm,
        enlargelimits=false,
        xlabel = {Inferred load [MW]},
        title = {$\alpha=10\%$},
        ylabel style = {yshift = -0.4cm},
        legend style={draw=none,fill=none,font=\scriptsize},
        legend cell align={left},
        ]
            \addplot[red,fill=red!30,hist,line width=0.025cm,opacity=1] table [y=x5]{results_plot_data/save_SDP_inference_results_global_sense.dat};
            \addplot[blue,fill=blue!30,hist,line width=0.025cm,opacity=0.75] table [y=x5]{results_plot_data/save_SDP_inference_results.dat};
            \addplot[mark=none,line width=0.035cm,color=red,dashed] plot coordinates {(9,0) (9,400)};
        \end{axis}
        
    \end{tikzpicture}
    
    \vspace{-0.2cm}
    \caption{SP-ADMM: Results of privacy attack on the load sited at node 20. The plots depict the distribution of inferred load by an adversary across 1000 simulation runs for different adjacency coefficient. The red and blue distributions are given for the global and local sensitivities $\Delta_{z}$ and $\overline{\delta}_{z}$, respectively.}
    \label{fig:SP_ADMM_attack}
\end{figure*}
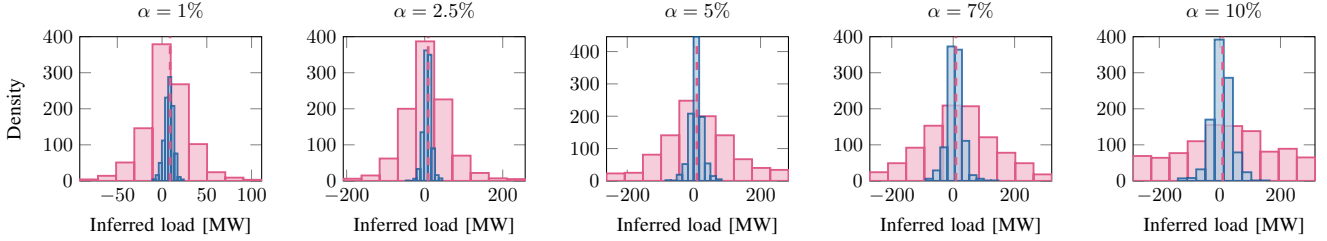

\pgfplotsset{
    compat=newest,
    /pgfplots/legend image code/.code={%
        \draw[mark repeat=2,mark phase=2,#1] 
            plot coordinates {
                (0cm,0cm) 
                (0.13cm,0cm)
                (0.16cm,0cm)
                (0.19cm,0cm)
                (0.3cm,0cm)%
            };
    },
}
\begin{figure}
\centering
\begin{tikzpicture}[scale=1, every node/.style={scale=0.75}]
    \begin{axis}[
        xshift = -4.75cm,
        width=4.5cm,
        height=4.5cm,
        colormap={bluewhite}{color=(blue!20) color=(red!50)},
        title = {$\xi_{z}\sim\text{Lap}(\frac{\delta_{z}^{k}}{\varepsilon})$},
        xticklabels={1,2,5,10,15},
        xtick={1,2,3,4,5},
        xtick style={draw=none},
        yticklabels={1.0,2.5,5.0,7.0,10.0},
        ylabel = {Adjacency $\alpha$, \%},
        ylabel style = {yshift = -0.3cm},
        xlabel = {Attack budget $T$},
        xlabel style = {yshift = 0.1cm},
        ytick={1,2,3,4,5},
        ymax=5.5,ymin=0.5,
        xmax=5.5,xmin=0.5,
        ytick style={draw=none},
        enlargelimits=false,
        point meta min=0.0,point meta max=16.7,
        nodes near coords={\pgfmathprintnumber\pgfplotspointmeta},
        nodes near coords style={
                yshift=-7pt,
        /pgf/number format/fixed,
                /pgf/number format/precision=1},
    ]
        \addplot[
            matrix plot,
            mesh/cols=5,
            point meta=explicit,draw=gray
        ] table [meta=C] {
            x           y           C
            % node 20 --- all information
        1.0000000	1.0000000	0.1600
        1.0000000	2.0000000	0.7300
        1.0000000	3.0000000	1.0500
        1.0000000	4.0000000	2.0600
        1.0000000	5.0000000	3.3000
        2.0000000	1.0000000	0.1500
        2.0000000	2.0000000	0.4965
        2.0000000	3.0000000	1.0183
        2.0000000	4.0000000	1.9274
        2.0000000	5.0000000	2.3474
        3.0000000	1.0000000	0.1392
        3.0000000	2.0000000	0.4457
        3.0000000	3.0000000	0.9764
        3.0000000	4.0000000	1.2747
        3.0000000	5.0000000	1.7552
        4.0000000	1.0000000	0.1223
        4.0000000	2.0000000	0.3864
        4.0000000	3.0000000	0.8096
        4.0000000	4.0000000	1.1616
        4.0000000	5.0000000	1.6645
        5.0000000	1.0000000	0.1141
        5.0000000	2.0000000	0.3864
        5.0000000	3.0000000	0.7768
        5.0000000	4.0000000	1.0646
        5.0000000	5.0000000	1.5214
        };
    \end{axis}

    \begin{axis}[
        width=4.5cm,
        height=4.5cm,
    colormap={bluewhite}{color=(blue!20) color=(red!50)},
        xticklabels={1,2,5,10,15},
        title = {$\xi_{z}\sim\text{Lap}(T\frac{\delta_{z}^{k}}{\varepsilon})$},
        xtick={1,2,3,4,5},
        xtick style={draw=none},
        yticklabels={1.0,2.5,5.0,7.0,10.0},
        ylabel = {Adjacency $\alpha$, \%},
        ylabel style = {yshift = -0.3cm},
        xlabel = {Attack budget $T$},
        xlabel style = {yshift = 0.1cm},
        ytick={1,2,3,4,5},
        ymax=5.5,ymin=0.5,
        xmax=5.5,xmin=0.5,
        ytick style={draw=none},
        enlargelimits=false,
        point meta min=0.0,point meta max=16.7,
        nodes near coords={\pgfmathprintnumber\pgfplotspointmeta},
        nodes near coords style={
                yshift=-7pt,
        /pgf/number format/fixed,
                /pgf/number format/precision=1},
    ]
        \addplot[
            matrix plot,
            mesh/cols=5,
            point meta=explicit,draw=gray
        ] table [meta=C] {
            x           y           C
            % node 20 --- all information
            1.0000000	1.0000000	0.1636192
            1.0000000	2.0000000	0.7397577
            1.0000000	3.0000000	1.0591781
            1.0000000	4.0000000	2.0659735
            1.0000000	5.0000000	3.3038126
            2.0000000	1.0000000	0.6228371
            2.0000000	2.0000000	1.6670010
            2.0000000	3.0000000	2.6217571
            2.0000000	4.0000000	4.3132903
            2.0000000	5.0000000	5.2873496
            3.0000000	1.0000000	1.1760945
            3.0000000	2.0000000	2.2619119
            3.0000000	3.0000000	3.8316252
            3.0000000	4.0000000	6.0297587
            3.0000000	5.0000000	8.4958609
            4.0000000	1.0000000	1.1807983
            4.0000000	2.0000000	2.3845995
            4.0000000	3.0000000	4.8000212
            4.0000000	4.0000000	7.5652106
            4.0000000	5.0000000	11.4075978
            5.0000000	1.0000000	1.0886097
            5.0000000	2.0000000	2.8598472
            5.0000000	3.0000000	6.5337714
            5.0000000	4.0000000	9.9006551
            5.0000000	5.0000000	16.5873872
        };
    \end{axis}
\end{tikzpicture}

\vspace{-0.25cm}
\caption{DP-ADMM: Mean absolute inference error, i.e., mismatch between the actual and inferred loads in MWh, across last $T$ iterations with (right) and without (left) application of Theorem \ref{th:compositon_two} for 100 simulation runs.}
\label{DP_ADMM_privacy_attack_results_matrix}
\end{figure}

The robustness of the algorithms to the load inference is assessed by
using the adversarial model in \eqref{OPF_attack_model}. The adversarial model identifies {\it all} network loads if the standard ADMM is used. The random perturbations of sub-problem solutions, however, prevent the adversary
from inferring the actual loads. The results focus on the load at bus
$20$, which has a median load in the first zone. As per Theorem
\ref{th:DP_privacy}, by specifying the adjacency coefficient $\alpha$, the
algorithms guarantee that, at a given iteration $k$, an adversary cannot
distinguish the magnitude of unknown load $d_{20}$ from any other
magnitude in the range $d_{20}\pm\alpha$. 

The load inference results for the DP-ADMM algorithm are shown in Fig. \ref{fig:DP_ADMM_attack}. The plots show the inferred load at every iteration of the algorithm assuming that only a single iteration is available to an adversary.  The inferred load is given as a probability density with each observation corresponding to a single iteration. By increasing $\alpha$, the inferred load deviates more substantially from the true value of 9 MW, hence, the probability of recovering the true load magnitude reduces. The load obfuscation with SP-ADMM is depicted in Fig. \ref{fig:SP_ADMM_attack}. Since the noise is fixed across iterations, the results display 1000 ADMM runs. Similarly to the DP-ADMM algorithm, increasing values of $\alpha$ result in wider distributions of inferred loads. It is important to note that the attacker observes only one sample from these distributions. Observe that the variance of load distributions is notably larger than that of DP-ADMM for a given adjacency value. Moreover, the support of the distributions in Fig. \ref{fig:SP_ADMM_attack} extends drastically with increasing values of $\alpha$, making the load inference essentially equivalent to a random guess. 

Fig. \ref{fig:SP_ADMM_attack} further shows that the use of upper bound on the global sensitivity $\Delta_{z}$ in SP-ADMM, which is set to the largest installed load in the system, results in much stronger privacy protection than the use of local sensitivity, which is set to be at least as much as the maximum local sensitivity observed across DP-ADMM iterations, i.e.,  $\overline{\delta}_{z}=\text{max}_{k}\{\delta_{z}^{k}\}_{k=1}^{K}$. Although both methods enable privacy protection, the formal privacy guarantees provided by SP-ADMM are only achieved with the use of global sensitivity. 

Finally, observe that every new iteration of DP-ADMM reveals more information to an adversary, as shown on the left plot in Fig. \ref{DP_ADMM_privacy_attack_results_matrix}. If the attack budget, i.e., the number of compromised iterations, increases up to $T$, an adversary recovers the load more precisely. To overcome this limitation, Theorem \ref{th:compositon_two} can be applied to preserve $\varepsilon-$differential privacy across $T$ iterations. This requires to scale the noise parameters by $T$. The corresponding results in the right plot in Fig. \ref{DP_ADMM_privacy_attack_results_matrix} show that the magnitude of the noise increases substantially, thus reducing the quality of load inference even with more information available to an adversary.

\subsubsection{Convergence analysis}  

The convergence statistics of the two algorithms obtained with 100
simulation runs are summarized in Fig. \ref{fig_convergence}. The
figure shows the evolution of the aggregated primal residual across
iterations highlighting important differences between dynamic and
static perturbations of sub-problem solutions. With dynamic
perturbations, the DP-ADMM algorithm perturbs the sub-problem solutions
at every iteration, and the magnitude of the noise increases with
$\alpha$. While the non-private ADMM converges in 59 iterations on
this test case, the DP-ADMM requires up to 300 iterations in average,
depending on the choice of $\alpha$. In contrast, the SP-ADMM
exhibits a similar computational complexity as the non-private ADMM
algorithm. Moreover, in average, the convergence of SP-ADMM is not
affected by the choice of adjacency coefficient.

\begin{figure}
\centering
% \resizebox{0.49\textwidth}{!}{%
\begin{tikzpicture}[font=\footnotesize]
\centering
\begin{axis}[
xshift=0cm,
small,
height=3.5cm,
width=5cm,
xmin=15,
xmax=80,
ymin=0,
ymax=15,
xlabel={iteration, $k$},
ylabel={$\norm{\tilde{\theta}^{k}-\overline{\theta}^{k+1}}_{2}$},
legend cell align=left,
legend pos= north east,
label style={font=\footnotesize},
every tick label/.append style={font=\scriptsize},
y label style={yshift=-0.2cm,font=\scriptsize},
legend style={yshift=-5,draw=none,fill=none,font=\scriptsize},
% ylabel shift = -5 pt,
% xlabel shift = -2.5 pt,
% xticklabels={0.2,0.4,0.6,0.8,1},xtick={200,400,600,800,1000},
% semithick,
tick style={major tick length=2pt,semithick,black},
separate axis lines,
axis x line*=bottom,
axis x line shift=5 pt,
axis y line*=left,
axis y line shift=5 pt,
]
% plot alpha 5 %
\addplot[name path=seven_one_5,draw=none,line width=0.005cm,color=blue!60] table [x index=0, y index=1] {results_plot_data/results_residual_0.05.dat};
\addplot[name path=seven_two_5,mark=none,line width=0.025cm,dashed,color=blue!60] table [x index=0, y index=2] {results_plot_data/results_residual_0.05.dat};
\addplot[name path=seven_three_5,mark=none,line width=0.005cm,color=blue!60] table [x index=0, y index=3] {results_plot_data/results_residual_0.05.dat};
\addplot[color=blue!60,opacity=0.5]fill between[of=seven_one_5 and seven_three_5];
% plot alpha 2.5 %
\addplot[name path=seven_one_25,draw=none,line width=0.005cm,color=blue!75] table [x index=0, y index=1] {results_plot_data/results_residual_0.025.dat};
\addplot[name path=seven_two_25,mark=none,line width=0.025cm,dashed,color=blue!75] table [x index=0, y index=2] {results_plot_data/results_residual_0.025.dat};
\addplot[name path=seven_three_25,mark=none,line width=0.005cm,color=blue!75] table [x index=0, y index=3] {results_plot_data/results_residual_0.025.dat};
\addplot[color=blue!75,opacity=0.5]fill between[of=seven_one_25 and seven_three_25];
% plot alpha 1 %
\addplot[name path=seven_one_1,mark=none,line width=0.005cm,color=blue!90] table [x index=0, y index=1] {results_plot_data/results_residual_0.01.dat};
\addplot[name path=seven_two_1,mark=none,line width=0.01cm,dashed,color=blue!90] table [x index=0, y index=2] {results_plot_data/results_residual_0.01.dat};
\addplot[name path=seven_three_1,mark=none,line width=0.005cm,color=blue!90] table [x index=0, y index=3] {results_plot_data/results_residual_0.01.dat};
\addplot[color=blue!90,opacity=0.5]fill between[of=seven_one_1 and seven_three_1];
\end{axis}
\begin{axis}[
xshift=4cm,
small,
height=3.5cm,
width=5cm,
xmin=10,
xmax=300,
ymin=0,
ymax=60,
xlabel={iteration, $k$},
ylabel={},
legend cell align=left,
legend pos= north east,
label style={font=\footnotesize},
every tick label/.append style={font=\scriptsize},
legend style={yshift=-5,draw=none,fill=none,font=\scriptsize},
% ylabel shift = -5 pt,
% xlabel shift = -2.5 pt,
% xticklabels={0.2,0.4,0.6,0.8,1},xtick={200,400,600,800,1000},
% semithick,
tick style={major tick length=2pt,semithick,black},
separate axis lines,
axis x line*=bottom,
axis x line shift=5 pt,
axis y line*=left,
axis y line shift=5 pt,
]
% plot alpha 10 %
\addplot[name path=seven_one,draw=none,line width=0.005cm,color=blue!30] table [x index=0, y index=1] {results_plot_data/results_residual_0.1.dat};
\addplot[name path=seven_two,mark=none,line width=0.025cm,dashed,color=blue!30] table [x index=0, y index=2] {results_plot_data/results_residual_0.1.dat};
\addplot[name path=seven_three,mark=none,line width=0.005cm,color=blue!30] table [x index=0, y index=3] {results_plot_data/results_residual_0.1.dat};
\addplot[color=blue!30,opacity=0.5]fill between[of=seven_one and seven_three];
% plot alpha 7 %
\addplot[name path=seven_one,mark=none,line width=0.005cm,color=blue!45] table [x index=0, y index=1] {results_plot_data/results_residual_0.07.dat};
\addplot[name path=seven_two,mark=none,line width=0.025cm,dashed,color=blue!45] table [x index=0, y index=2] {results_plot_data/results_residual_0.07.dat};
\addplot[name path=seven_three,mark=none,line width=0.005cm,color=blue!45] table [x index=0, y index=3] {results_plot_data/results_residual_0.07.dat};
\addplot[color=blue!45,opacity=0.5]fill between[of=seven_one and seven_three];
\end{axis}
\centering
\begin{axis}[
xshift=0cm,
yshift=-3cm,
small,
height=3.5cm,
width=9cm,
xmin=15,
xmax=63,
ymin=0,
ymax=10,
xlabel={iteration, $k$},
ylabel={$\norm{\tilde{\theta}^{k}-\overline{\theta}^{k+1}}_{2}$},
legend cell align=left,
legend pos= north east,
label style={font=\footnotesize},
every tick label/.append style={font=\scriptsize},
y label style={yshift=-0.2cm,font=\scriptsize},
legend columns=-1,
legend style={yshift=-3cm,xshift=0.3cm,fill=none,font=\scriptsize},
tick style={major tick length=2pt,semithick,black},
separate axis lines,
axis x line*=bottom,
axis x line shift=5 pt,
axis y line*=left,
axis y line shift=5 pt,
]

% plot alpha 10 %
\addplot[name path=seven_one,draw=none,line width=0.015cm,color=blue!30,forget plot] table [x index=0, y index=1] {results_plot_data/results_SPD_residual_0.1.dat};
\addplot[name path=seven_three,mark=none,line width=0.015cm,color=blue!30,forget plot] table [x index=0, y index=3] {results_plot_data/results_SPD_residual_0.1.dat};
\addplot[color=blue!30]fill between[of=seven_one and seven_three];\addlegendentry{$\alpha=10$};
% plot alpha 7 %
\addplot[name path=seven_one,mark=none,line width=0.015cm,color=blue!45,forget plot] table [x index=0, y index=1] {results_plot_data/results_SPD_residual_0.07.dat};
\addplot[name path=seven_three,mark=none,line width=0.015cm,color=blue!45,forget plot] table [x index=0, y index=3] {results_plot_data/results_SPD_residual_0.07.dat};
\addplot[color=blue!45]fill between[of=seven_one and seven_three];\addlegendentry{$\alpha=7$};
% plot alpha 5 %
\addplot[name path=seven_one_5,line width=0.015cm,color=blue!60,forget plot] table [x index=0, y index=1] {results_plot_data/results_SPD_residual_0.05.dat};
\addplot[name path=seven_three_5,mark=none,line width=0.015cm,color=blue!60,forget plot] table [x index=0, y index=3] {results_plot_data/results_SPD_residual_0.05.dat};
\addplot[color=blue!60,opacity=1]fill between[of=seven_one_5 and seven_three_5];\addlegendentry{$\alpha=5$};
% plot alpha 2.5 %
\addplot[name path=seven_one_25,line width=0.015cm,color=blue!75,forget plot] table [x index=0, y index=1] {results_plot_data/results_SPD_residual_0.025.dat};
\addplot[name path=seven_two_25,mark=none,line width=0.025cm,dashed,color=blue!75,forget plot] table [x index=0, y index=2] {results_plot_data/results_SPD_residual_0.025.dat};
\addplot[name path=seven_three_25,mark=none,line width=0.015cm,color=blue!75,forget plot] table [x index=0, y index=3] {results_plot_data/results_SPD_residual_0.025.dat};
\addplot[color=blue!75,opacity=1]fill between[of=seven_one_25 and seven_three_25];\addlegendentry{$\alpha=2.5$};
% plot alpha 1 %
% \addplot[name path=seven_one_1,mark=none,line width=0.015cm,color=blue!90,forget plot] table [x index=0, y index=1] {results_plot_data/results_SPD_residual_0.01.dat};
\addplot[name path=seven_two_1,mark=none,line width=0.025cm,dashed,color=blue!90,forget plot] table [x index=0, y index=2] {results_plot_data/results_SPD_residual_0.01.dat};
% \addplot[name path=seven_three_1,mark=none,line width=0.015cm,color=blue!90,forget plot] table [x index=0, y index=3] {results_plot_data/results_SPD_residual_0.01.dat};
\addplot[color=blue!90,opacity=1]fill between[of=seven_one_1 and seven_three_1];\addlegendentry{$\alpha=1$};

\end{axis}
\end{tikzpicture}
% }

\vspace{2.75cm}
\caption{Convergence of DP-ADMM (top) and SP-ADMM (bottom) algorithms on the 3-zone IEEE 118-node system for different adjacency coefficient $\alpha$ in \%. The dashed lines indicate the average residual across 100 runs, whereas the colored areas indicate the spread between the minimum and maximum values of the residual at iteration $k$. Best view in colors.}
\label{fig_convergence}
\end{figure}
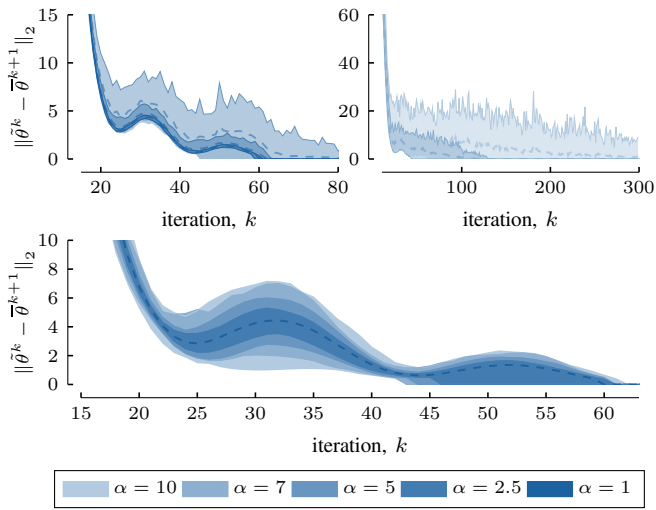

\subsubsection{Fidelity analysis}
It remains to quantify the loss in efficiency of
differentially-private OPF solutions. The average optimality loss induced by DP-ADMM and SP-ADMM algorithms for 100 runs is provided in Table \ref{fig_convergence}. The results for attack budget $T=1$ show that with increasing privacy requirements, both algorithms converge in sub-optimal solutions as compared to the non-private ADMM solution. However, due to dynamically updated zero-mean perturbations, the DP-ADMM has notability better fidelity than the SP-ADMM, which fixes the noise across iterations and constantly steers the OPF dispatch from the optimal solution. This unfolds the following trade-offs among the algorithms: despite better convergence properties of SP-ADMM, it yields a larger optimality gap compared to the DP-ADMM, which demonstrates weaker convergence statistics.

\begin{table}[]
\caption{Optimality loss induced by DP-ADMM and SP-ADMM $(\%)$}
\centering
\begin{tabular}{l|ccccc}
\specialrule{0.70pt}{0.70pt}{0.70pt}
$\alpha,$\%   & 1    & 2.5  & 5    & 7     & 10 \\
\specialrule{0.70pt}{0.70pt}{0.70pt}
DP-ADMM & 0.48 & 0.92 & 1.23 & 1.51  & 3.83  \\
SP-ADMM & 0.28 & 4.33 & 11.0 & 11.35 & 20.41 \\
\specialrule{0.70pt}{0.70pt}{0.70pt}
\end{tabular}
\end{table}

%%%%%%%%%%%%%%%%%%%%%%%%%%%%%%%%%%%%%%%%%%%%%%%%%%%%%%%%
% \fi
%%%%%%%%%%%%%%%%%%%%%%%%%%%%%%%%%%%%%%%%%%%%%%%%%%%%%%%%

\section{Conclusions}
Although the distributed algorithms have been long trusted to preserve the privacy of network parameters in OPF computations, this paper shows that the standard distributed algorithms do not ensure the information integrity as the sensitive parameters, e.g., electrical loads, are leaked through the exchange of coordination signals. To overcome this limitation, this paper introduces two privacy-preserving OPF ADMM algorithms that satisfy the definition of $\varepsilon$-differential privacy. The algorithms provide privacy be means of either static or dynamic perturbations of the sub-problem solutions at each iteration. The paper shows theoretically and through numerical results that the two algorithms are able to negate the adversarial inference of sensitive information from coordination signals. Despite their complementary privacy properties, the numerical performance of the two algorithms is mutually exclusive: if static perturbations demonstrate a more robust convergence, their fidelity with respect to the non-private solution is lower than that of dynamic perturbations with weaker convergence statistics.

\bibliographystyle{IEEEtran}
\bibliography{references}
% \balance
\endgroup
\end{document}